\documentclass[ejs]{imsart}
\RequirePackage[OT1]{fontenc}
\RequirePackage{amsthm,amsmath,natbib}
\RequirePackage[colorlinks]{hyperref}
\RequirePackage{hypernat}
\usepackage{latexsym}
\usepackage{amssymb}
\usepackage[dvips,usenames]{color}
\usepackage[dvips]{graphicx}
\usepackage{multicol}

\pubyear{2008}

\startlocaldefs
\newtheorem{teo}{Theorem}[section]

\newtheorem{lem}[teo]{Lemma}

\newcommand{\bfun}{\left\{\begin{array}{ll}}
\newcommand{\efun}{\end{array}\right.}
\newcommand{\bfunB}{\begin{array}{ll}}
\newcommand{\efunB}{\end{array}}
\newcommand{\Prob}{\textsf{P}}
\newcommand{\card}{\textsf{card}}

\newcommand{\MyR}{\mathbb{R}}
\newcommand{\MyZ}{\mathbb{Z}}

\newcommand{\MyNN}{\mathbb{N}}

\newcommand{\MyI}{\textsf{I}}
\newcommand{\nSeq}{\{1,\ldots,n\}}

\newcommand{\Exp}{\mathbb{E}}

\newcommand{\bx}{\mathbf{x}}
\newcommand{\bX}{\mathbf{X}}
\newcommand{\bz}{\mathbf{z}}

\newcommand{\Beso}{\mathcal{B}}
\newcommand{\Lp}{\textsf{L}}
\newcommand{\Ldue}{\textsf{L}^2}

\newcommand{\MyJ}{\mathrm{J}}

\renewcommand{\epsilon}{\varepsilon}

\newcommand{\bi}{\begin{itemize}}
\newcommand{\ei}{\end{itemize}}
\newcommand{\be}{\begin{eqnarray*}}
\newcommand{\ee}{\end{eqnarray*}}
\newcommand{\bequ}{\begin{equation}}
\newcommand{\eequ}{\end{equation}}
\endlocaldefs

\begin{document}

\begin{frontmatter}
\title{Warped Wavelet and Vertical Thresholding}
\runtitle{Warped Wavelet and Vertical Thresholding}

\begin{aug}
\author{\fnms{Pierpaolo} \snm{Brutti}\ead[label=e1]{pbrutti@luiss.it}
\ead[label=u1,url]{http://docenti.luiss.it/brutti}
}

\address{Dipartimento di Scienze Economiche e Aziendali\\
LUISS, Guido Carli \\
Viale Romania 32, 00197 Roma, Italy. \\
\printead{e1}\\
\printead{u1}
}

\runauthor{P. Brutti}

\affiliation{LUISS, Guido Carli}

\end{aug}

\begin{abstract}
Let $\{(X_i,Y_i)\}_{i\in\nSeq}$ be an i.i.d. sample from the random design regression model $Y = f(X) + \epsilon$ with $(X,Y) \in [0,1] \times [-M,M]$. In dealing with such a model, adaptation is naturally to be intended in terms of $\Lp^2([0,1],G_X)$ norm where $G_X(\cdot)$ denotes the (known) marginal distribution of the design variable $X$. Recently much work has been devoted to the construction of estimators that adapts in this setting (see, for example,
\cite{Binev:Cohen:Dahmen:DeVore:Temlyakov:2004,DeVore:Kerk:Picard:Temlyakov:2004b,DeVore:Kerk:Picard:Temlyakov:2004a,Gyorfi:etal:2002}),
but only a few of them come along with a easy--to--implement computational scheme.
Here we propose a family of estimators based on the \emph{warped wavelet} basis recently introduced by Picard and Kerkyacharian \cite{KP04} and a tree-like thresholding rule that takes into account the hierarchical (across-scale) structure of the wavelet coefficients.
We show that, if the regression function belongs to a certain class of approximation spaces defined in terms of $G_X(\cdot)$, then our procedure is adaptive and converge to the true regression function with an optimal rate. The results are stated in terms of excess probabilities as in \cite{Cucker:Smale:2002}.
\end{abstract}

\begin{keyword}[class=AMS]
\kwd[Primary ]{62G07}
\kwd{60K35}
\kwd[; secondary ]{62G20}
\end{keyword}

\begin{keyword}
\kwd{Regression with random design}
\kwd{Wavelets}
\kwd{Block thresholding}
\kwd{Warped Wavelets}
\kwd{Adaptive Approximation}
\kwd{Universal Algorithms}
\kwd{Muckenhoupt weights}
\end{keyword}
\end{frontmatter}

\section{Introduction}\label{sec:CH3:WarpedWavelet}

Wavelet bases are ubiquitous in modern nonparametric statistics starting from the 1994 seminal paper by Donoho and Johnstone \cite{Dono:John:1994}. What makes them so appealing to statisticians is their ability to capture the relevant features of smooth signals in a few ``big'' coefficients at high scales (low frequencies) so that zero thresholding the small ones, results in an effective denoising scheme (see \cite{Vida:1999}).

Although these well known results about thresholding techniques were usually obtained assuming a fixed (and possibly equispaced) design \cite{Dono:John:1994,Doh:John:Kerk:Picard:1995}, it was quite reassuring to see how they carry over almost unchanged to the irregular design case. As a matter of fact, in the case of irregular design, various attempts to solve this problem has been made: see, for instance, the interpolation methods of Hall and Turlach \cite{Hall:Turlach:1997} and Kovac and Silverman \cite{Kovac:Silverman:2000}; the binning method of Antoniadis et al. \cite{Antoniadis:etal:1997}; the transformation method of Cai and Brown \cite{Cai:Brown:1998}, or its recent refinements by Maxim \cite{Maxim:2002} for a random design; the weighted wavelet transform of Foster \cite{Foster:1996}; the isometric method of Sardy et al. \cite{Sardy:etal:1999}; the penalization method of Antoniadis and Fan \cite{Anto:Fan:2001}; and the specific construction of wavelets adapted to the design of Delouille et al. \cite{Deloui:Franke:vonSa:2001,Deloui:Simo:vonSa:2004} and Jansen et al. \cite{Vanr:Jans:Bult:2002}. See also Pensky and Vidakovic \cite{Pensky:Vida:2001}, and the monograph \cite{Gyorfi:etal:2002}.

The main drawback common to most of the methods just mentioned can be found, with no surprise, on the computational side: compared, for instance, with the usual thresholding technique, the calculations are, in general, less direct. To fix this problem, Kerkyacharian and Picard \cite{KP04} propose \emph{warped} wavelet basis. The idea is as follow. For a signal observed at some design points, $Y(t_i)$, $i \in \{1,\ldots,2^\mathrm{J}\}$, if the design is regular ($t_k = k/2^\mathrm{J}$), the standard wavelet decomposition algorithm starts with $s_{\mathrm{J},k} = 2^{\mathrm{J}/2}Y(k/2^{\mathrm{J}})$ which approximates the scaling coefficient $\int Y(x)\phi_{\mathrm{J},k}(x)\mathrm{d}x$, with $\phi_{\mathrm{J},k}(x) = 2^{\mathrm{J}/2}\phi(2^\mathrm{J} x - k)$ and $\phi(\cdot)$ the so--called scaling function or father wavelet (see \cite{Mallat:1998} for further information). Then the cascade algorithm is employed to obtain the wavelet coefficients $d_{j,k}$ for $j \leqslant \mathrm{J}$, which in turn are thresholded. If the design is not regular, and we still employ the \emph{same} algorithm, then for a function $H(\cdot)$ such that $H(k/2^\mathrm{J}) = t_k$, we have $s_{\mathrm{J},k} = 2^{\mathrm{J}/2}Y(H(k/2^{\mathrm{J}}))$. Essentially what we are doing is to decompose, with respect to a standard wavelet basis, the function $Y(H(x))$ or, if $G \circ H(x) \equiv x$, the original function $Y(x)$ itself but with respect to a new \emph{warped} basis $\{\psi_{j,k}(G(\cdot))\}_{(j,k)}$.
In the regression setting, this means replacing the standard wavelet expansion of the function $f(\cdot)$ by its expansion on the new basis $\{\psi_{j,k}(G(\cdot))\}_{(j,k)}$, where $G(\cdot)$ is adapting to the design: it may be the distribution function of the design, or its estimation, when it is unknown (not our case). An appealing feature of this method is that it does not need a new algorithm to be implemented: just standard and widespread tools.
Of course the properties of this basis depend on the warping factor $G(\cdot)$. In \cite{KP04} the authors provide the conditions under which this new basis behaves, at least for statistical purposes, as well as ordinary wavelet bases with respect to $\Lp^p([0,1],\mathrm{d}x)$ norms with $p\in(0,+\infty)$. This condition properly quantifies the departure from the uniform distribution and happens to be associated with the notion of Muckenhoupt weights (see \cite{Garcia:Martell:2001,Stein:1993}).

Now the problem is that we do not need good estimators in $\Lp^p([0,1],\mathrm{d}x)$. What we need are (easy to compute) estimators that adapt in $\Lp^2([0,1],G_X)$. As a matter of fact it is possible to prove that the main results contained in \cite{KP04} can be extended to this new setting once we assume $G_X(\cdot)$ to be known as in \cite{Chesneau:2007}, the case of an unknown $G_X(\cdot)$ being beyond the scope of this work (see \cite{Kerk:Picard:2005}).

Here we propose a particular variation on the basic thresholding procedure advanced in \cite{KP04}, that can be motivated as follow.
\begin{figure}[htb]
\begin{center}
\includegraphics[keepaspectratio=true,width=10cm,clip=true]{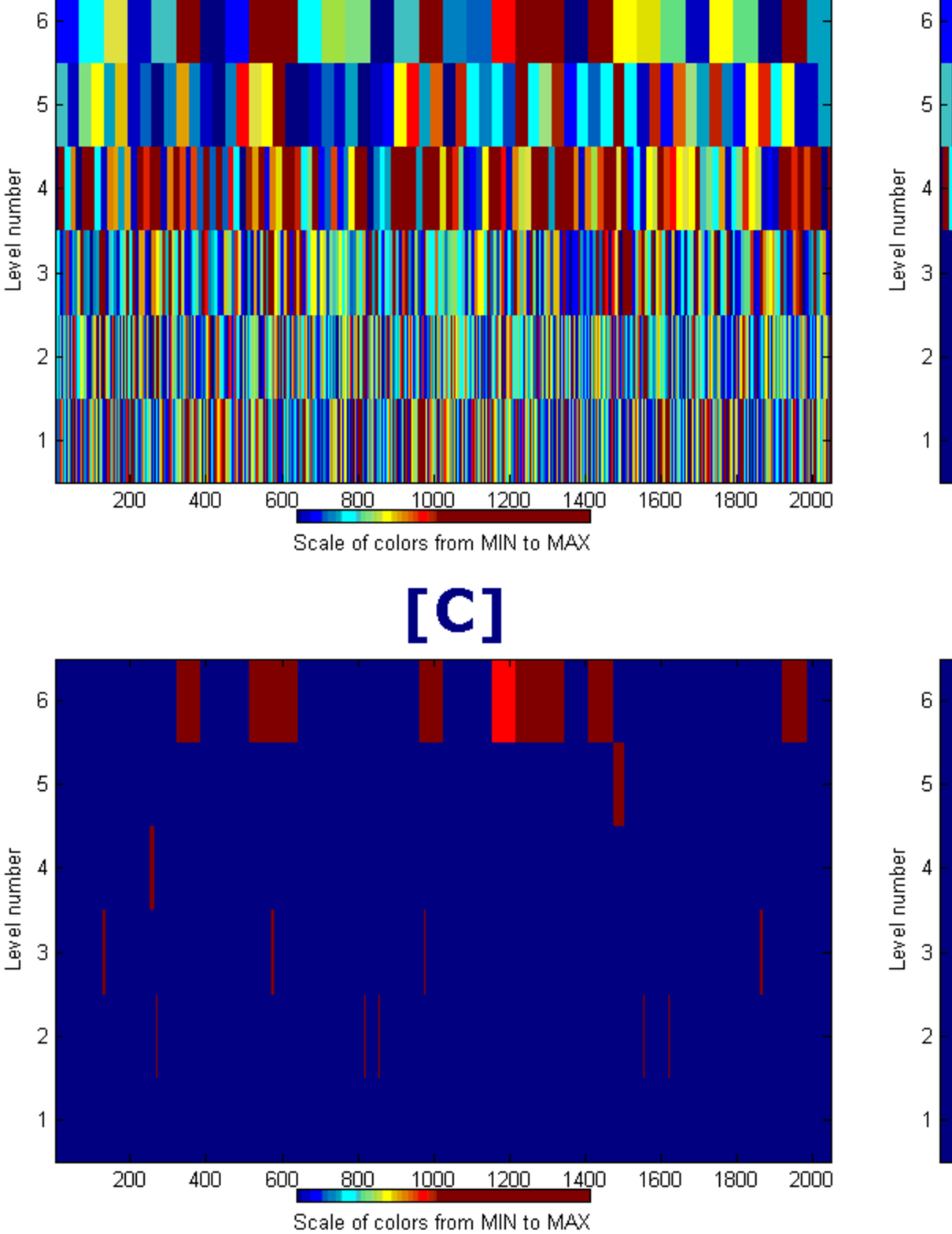}
\end{center}
\caption[Examples of thresholding rules]{\textsf { \small Examples of thresholding rules: [A] - Original wavelet coefficients; [B] - Linear thresholding; [C] - Nonlinear (hard) thresholding; [D] - Vertical (hard) thresholding.
}}
\label{fig:CH3:figura1}
\end{figure}
In a variety of real--life signals, significant wavelet coefficients often occur in clusters at adjacent scales and locations. Irregularities, like a discontinuity for example, in general tend to affect the whole block of coefficients corresponding to wavelet functions whose ``support'' contains them. For this reason it is reasonable to expect that the risk of ``blocked'' thresholding rules might compare quite favorably with other classical estimators based on level--wise or global thresholds. The literature is filled with successful examples of ``horizontally'' (within scales) blocked rules derived from both, purely frequentist arguments \cite{Cai:1999,Cai:Silver:2001,Chesneau:2007,Hall:Kerk:Picard:1999}, or Bayesian reasonings of some flavor \cite{Abra:Besb:Sapa:2002,DeCand:Vida:2001,Wang:Wood:2003}. Recently, an increasing amount of work has been devoted to study a new class of ``vertically'' (across scales, see Figure \ref{fig:CH3:figura1}) blocked or \emph{treed} rules \cite{Autin:Picard:Riv:2004,Io:2004,Cohen:Dah:Daub:deVore:2001,Fryz:2004,Rom:Choi:Bara:2001}, that have proved to be of invaluable help in at least two settings of great importance: the construction of \emph{adaptive} pointwise confidence intervals \cite{Picard:Trib:2000} and the derivation of pointwise estimators of a regression function that adapt rate optimally under what we could call a \emph{focused} performance measure \cite{Cai:Low:2005a}.

For this reason, adapting some techniques developed in \cite{Binev:Cohen:Dahmen:DeVore:Temlyakov:2004} to the current (simplified) setting, in Section \ref{subsec:CH3:TreeApproximation} we show how vertically zero--thresholding the \emph{warped} wavelet coefficients actually results in an universal smoother with good properties in $\Lp^2([0,1],G_X)$ over reasonably large approximation spaces.

\section{Tree--Structured Warped Approximations}\label{subsec:CH3:TreeApproximation}

We shall now discuss in greater details nonlinear approximation processes based on warped
wavelet bases where a tree structure is pre--imposed on the preserved coefficients.
We will start following closely \cite{Binev:Cohen:Dahmen:DeVore:Temlyakov:2004} by reviewing some basic facts about partitions and how they are related to adaptive approximation.
Then we present the universal algorithm based on adaptive partitions coming from a warped wavelet decomposition and its theoretical properties.

In the spirit of the recent paper by Cucker and Smale \cite{CS}, we will measure the performances of our estimator by studying its convergence both in probability and expectation. More specifically, let $\Prob\{\cdot\}$ be a -- generally unknown or partially unknown-- Borel measure defined on $\mathcal{Z} = \mathcal{X} \times \mathcal{Y} \subset \MyR^d \times \MyR$, and consider again a nonparametric regression problem where we want to estimate the conditional mean $f(\bx) = \Exp(Y|\bX = \bx)$ from an i.i.d. sample of size $n$, $\bz = \mathbf{z}_n = \{(\bx_i,y_i)\}_{i\in\nSeq}$, drawn from $\Prob\{\cdot\}$. Assume further that, chosen an \emph{hypothesis space} $\mathcal{H}$ from which our candidate estimators $f_{\mathbf{z}}(\cdot)$ comes from, we shall measure the approximation error of $f_{\mathbf{z}}(\cdot)$ in the $\Lp^2(\mathcal{X},G_\bX)$ norm, where $G_\bX(\cdot)$ is the (marginal) distribution of the design variable $\bX$. Here, as in the previous section, we will assume $G_\bX(\cdot)$ to be known. So, given $f_{\mathbf{z}} \in \mathcal{H}$, the quality of its performance is measured by
\[
\|f - f_{\bz}\| = \|f - f_{\bz} \|_{\Lp^2(\mathcal{X},G_\bX)}.
\]
Clearly this quantity is stochastic in nature and, consequently, it is generally not possible to say anything about it for a fixed $\mathbf{z}$. Instead we look at the behavior in
\emph{probability} as measured by
\[
\Prob^{ \otimes n} \left\{ {\mathbf{z}: \| {f - f_{\mathbf{z}} } \| > \eta } \right\},\;\eta > 0
\]
or the expected error
\[
\Exp^{\otimes n} \big( { \| {f - f_{\bz} } \|} \big) =
    \int { \| {f - f_{\bz} } \|\mathrm{d}\Prob^{ \otimes n} },
\]
where $\Prob^{\otimes n}\{\cdot\}$ denotes the $n$--fold tensor product of $\Prob\{\cdot\}$. Clearly, given a bound for $\Prob^{ \otimes n} \left\{ {\mathbf{z}: \| {f - f_{\mathbf{z}} } \| > \eta } \right\}$, we can immediately obtain another bound for the expected error since
\bequ\label{eq:CH3:Exp2Prob}
\Exp^{ \otimes n} \big( {\| {f - f_\bz }\|} \big) = \int_0^{ + \infty } {\Prob^{ \otimes n} \left\{ {\bz: \| {f - f_{\bz} } \| > \eta } \right\}\mathrm{d}\eta }.
\eequ
As we will see in Section \ref{sec:CH3:ProofsWarped}, bounding probabilities like $\Prob^{\otimes n}\{\cdot\}$ usually requires some kind of concentration of measure inequalities (see \cite{Boucheron:Bousquet:Lugosi:2004:A}).

Now, suppose that we have chosen a reasonable hypothesis space $\mathcal{H}$. We still need to address the problem of how to find an estimator $f_\bz(\cdot)$ for the regression function $f(\cdot)$. One of the most widespread criteria (see \cite{CS,DeVore:Kerk:Picard:Temlyakov:2004b,Gyorfi:etal:2002}, and references therein) is the so called \emph{empirical risk minimization} (least--square data fitting).

Empirical risk minimization is motivated by the fact that the regression function $f(\cdot)$ is the minimizer of
\[
\mathcal{E}(w) = \int {\left[ {w(\bx) - y} \right]^2 \mathrm{d}\Prob}.
\]
That is
\[
\mathcal{E}(f) = \mathop {\inf }\limits_{w \in \Ldue(\mathcal{X},G_\bX )} \mathcal{E}(w).
\]
This suggests to consider the problem of minimizing the empirical loss
\[
\mathcal{E}_\bz (w) = \frac{1}{n}\sum\limits_{i = 1}^n {[w(\bx_i ) - y_i ]^2 },
\]
over all $w \in \mathcal{H}$. So, in the end, we found an implementable form for our candidate estimator
\[
f_\bz = f_{\bz,\mathcal{H}}=\mathop {\arg \min }\limits_{w \in \mathcal{H}}\mathcal{E}_\bz(w),
\]
the so--called \emph{empirical minimizer}. Notice that given a finite ball in a linear or nonlinear finite dimensional space, the problem of finding $f_\bz(\cdot)$ is numerically solvable.

In the following we will see how to build the hypothesis space $\mathcal{H}$ from refinable partitions of the design space $\mathcal{X}$ and then, how this is related to (warped) wavelet basis. Typically $\mathcal{H} = \mathcal{H}_n$ depends on a finite number $\MyJ(n)$ of parameters as, for example, the dimension of a linear space or, equivalently, the number of basis functions we use to generate it. In many cases, this number $\MyJ$ is chosen using some a priori assumption on the regression function. In other procedures, the number $\MyJ$ avoids any a priori assumptions by adapting to the data. We shall be interested in estimators of the latter type.

\subsection{Partitions, Adaptive Approximation and Least--Squares Fitting}\label{subsubsec:CH3:PartitionsAndAdapProxy}

We will now review some basic facts about partitions and how they are related to adaptive approximation. The treatment follows closely \cite{Binev:Cohen:Dahmen:DeVore:Temlyakov:2004}. A partitions $\Lambda$ of $\mathcal{X} \subset [0,1]^d$ is usually built through a refinement strategy. We first describe the prototypical example of dyadic partitions and then, in the following section, we will make the link with orthonormal expansions through a wavelet basis. So let $\mathcal{X} = [0,1]^d$, and denote by $\mathcal{D}_j = \mathcal{D}_j(\mathcal{X})$ the collection of dyadic subcubes of $\mathcal{X}$ of sidelength $2^{-j}$ and $\mathcal{D} = \bigcup\nolimits_{j = 0}^\infty {\mathcal{D}_j}$. These cubes are naturally aligned on a tree $\mathcal{T} = \mathcal{T}(\mathcal{D})$. Each node of the tree $\mathcal{T}$ is a cube $\textsf{I} \in \mathcal{D}$. If $\textsf{I} \in \mathcal{D}_j$, then its children are the $2^d$ dyadic cubes of $\textsf{J} \in \mathcal{D}_{j+1}$ with $\textsf{J} \subset \textsf{I}$. We denote the set of children of \textsf{I} by $\mathcal{C}(\textsf{I})$. We call \textsf{I} the parent of each such child \textsf{J} and write $\textsf{I} = \mathcal{P}(\textsf{J})$. The cubes in $\mathcal{D}_j(\mathcal{X})$ form a uniform partition in which every cube has the same measure $2^{-j\,d}$.

More in general, we say that a collection of nodes $\widetilde{\mathcal{T}}$ is a \emph{proper} subtree of $\mathcal{T}$ if:
\begin{itemize}
\item the root node $\textsf{I} \equiv \mathcal{X}$ is in $\widetilde{\mathcal{T}}$,
\item if $\textsf{I} \neq \mathcal{X}$ is in $\widetilde{\mathcal{T}}$ then its parent $\mathcal{P}(\textsf{I})$ is also in $\widetilde{\mathcal{T}}$.
\end{itemize}
Any finite proper subtree $\widetilde{\mathcal{T}}$ is associated to a unique partition $\Lambda = \Lambda(\widetilde{\mathcal{T}})$ which consists of its outer leaves, by which we mean those $\textsf{J} \in \mathcal{T}$ such that $\textsf{J} \notin \widetilde{\mathcal{T}}$ but $\mathcal{P}(\textsf{J})$ is in $\widetilde{\mathcal{T}}$. One way of generating adaptive partitions is through some refinement strategy. One begins at the root $\mathcal{X}$ and decides whether to refine $\mathcal{X}$ (i.e. subdivide $\mathcal{X}$) based on some refinement criteria. If $\mathcal{X}$ is subdivided, then one examines each child and decides whether or not to refine such a child based on the refinement strategy.

We could also consider more general refinements. Assume, for instance, that $a \geqslant 2$ is a fixed integer. We assume that if $\mathcal{X}$ is to be refined, then its children consist of $a$ subsets of $\mathcal{X}$ which are a partition of $\mathcal{X}$. Similarly, for each such child there is a rule which spells out how this child is refined. We assume that the child is also refined into $a$ sets which form a partition of the child. Such a refinement strategy also results in a tree $\mathcal{T}$ (called the \emph{master tree}) and children, parents, proper trees and partitions are defined as above for the special case of dyadic partitions. The refinement level $j$ of a node is the smallest number of refinements (starting at root) to create this node. Note that to describe these more general refinements in terms of basis functions, we need to introduce the concept of \emph{warped} multi--wavelets and wavelet packets, but this is beyond the scope of the present work.

We denote by $\mathcal{T}_j$ the proper subtree consisting of all nodes with level $< j$ and we denote by $\Lambda_j$ the partition associated to $\mathcal{T}_j$, which coincides with
$\mathcal{D}_j(\mathcal{X})$ in the above described dyadic partition case. Note that in contrast to this case, the $a$ children may not be similar in which case the partitions $\Lambda_j$ are not spatially uniform (we could also work with even in more generality and allow the number of children to depend on the cell to be refined, while remaining globally bounded by some fixed $a$). It is important to note that the cardinalities of a proper tree $\widetilde{\mathcal{T}}$ and of its associated partition $\Lambda(\widetilde{\mathcal{T}})$ are equivalent. In fact one easily checks that
\[
\textsf{card}\big( \Lambda(\widetilde{\mathcal{T}}) \big) = (a - 1) \, \textsf{card}\big( \widetilde{\mathcal{T}}\big) + 1,
\]
by remarking that each time a new node gets refined in the process of building an adaptive partition, $\textsf{card}(\widetilde{\mathcal{T}})$ is incremented by $1$ and $\textsf{card}(\Lambda)$ by $a-1$.

Given a partition $\Lambda$, we can easily use it to approximate functions supported on $\mathcal{X}$. More specifically, let us denote by $\mathcal{S}_\Lambda$ the space of piecewise constant functions -- normalized in $\Ldue(\mathcal{X},G_X)$ -- subordinate to $\Lambda$. Each $f \in \mathcal{S}_\Lambda$ can then be written as
\[
f( \cdot ) = \sum\limits_{\textsf{I} \in \Lambda } {c_\textsf{I} \tfrac{1}{{\sqrt {G_X (\textsf{I})} }}\boldsymbol{1}_\textsf{I} ( \cdot )},
\]
where $\boldsymbol{1}_\textsf{I}(\cdot)$ denotes the indicator function of any set $\textsf{I} \subset \mathcal{X}$. The best approximation of a given function $f \in \Ldue(\mathcal{X},G_X)$ by the elements of $\mathcal{S}_\Lambda$ is given by
\[
\Pi _\Lambda(f)( \cdot ) = \sum\limits_{\textsf{I} \in \Lambda } {s_\textsf{I}\tfrac{1}{{\sqrt {G_X (\textsf{I})} }} \boldsymbol{1}_\textsf{I} ( \cdot )},
\]
where
\bequ\label{eq:CH3:piecewiseConsts_coeff}
s_\textsf{I}  = \left\langle {f,\tfrac{1}{{\sqrt {G_X (\textsf{I})} }}\boldsymbol{1}_\textsf{I} } \right\rangle _{\Ldue(G_X )},
\eequ
and $s_\textsf{I} \equiv 0$ in case $G_X(\textsf{I}) \equiv 0$.

In practice, we can consider two types of approximations corresponding to \emph{uniform refinement} and
\emph{adaptive refinement}. We first discuss uniform refinement. Let
\[
\mathcal{E}_\MyJ(f) = \|f - \Pi_{\Lambda_\MyJ}(f)\|_{\Ldue(G_\bX)}, \quad \MyJ \in \MyNN_0,
\]
which is the error for uniform refinement. The decay of this error to zero is connected with the smoothness of $f(\cdot)$ as measured in $\Ldue(\mathcal{X},G_\bX)$. We shall denote by $\mathcal{A}^s$ the approximation space (see the review in \cite{DeVore:1998}), consisting of all functions $f \in \Ldue(\mathcal{X},G_\bX)$ such that
\bequ\label{eq:CH3:binev_(3)}
\mathcal{E}_\MyJ(f) \leqslant M_0 a^{-\MyJ\,s}, \quad \MyJ \in \MyNN_0.
\eequ
Notice that $\card(\Lambda_\MyJ) = a^{\MyJ}$, so that the decay in Equation (\ref{eq:CH3:binev_(3)}) is like $N^{-s}$ with $N$ the number of elements in the
partition. The smallest $M_0$ for which Equation (\ref{eq:CH3:binev_(3)}) holds serves to define the semi-norm $| f |_{\mathcal{A}^s}$ on $\mathcal{A}^s$. The space $\mathcal{A}^s$ can be viewed as a smoothness space of order $s > 0$ with smoothness measured with
respect to $G_\bX(\cdot)$.
For example, if $G_\bX(\cdot)$ is the Lebesgue measure and we use dyadic partitioning then $\mathcal{A}^{s/d} = \Beso^{2,s}_\infty$, $ s \in (0,1]$, with equivalent norms. Here $\Beso^{2,s}_\infty$ is the Besov space which can be described in terms of the differences as
\[
\|w(\cdot + h) - w(\cdot)\|_{\Ldue(\mathrm{d}x)} \leqslant M_0 |h|^s, \quad x,h\in\mathcal{X}.
\]

Instead of working with a--priori fixed partitions there is a second kind of approximation where the partition is generated adaptively and will vary with $f(\cdot)$ . Adaptive partitions are typically generated by using some refinement criterion that determines whether or not to subdivide a given cell. We shall consider a refinement criteria that was introduced to build adaptive wavelet constructions such as those given by Cohen \emph{et al.} in \cite{Cohen:Dah:Daub:deVore:2001} for image compression. This criteria  is analogous to thresholding wavelet coefficients. Indeed, it would be exactly this criteria if we were
to construct a wavelet (\emph{Haar} like) bases for $\Ldue(\mathcal{X},G_\bX)$.
For each cell \textsf{I} in the master tree $\mathcal{T}$ and any $w\in\Ldue(\mathcal{X},G_\bX)$ we define
\bequ\label{eq:CH3:Theoretical_residuals}
\nu_\textsf{I} = \nu_\textsf{I} (w) = \sqrt { {\sum\limits_{\textsf{J} \in \mathcal{C}(\textsf{I})} {s_\textsf{J}^2 } } - s_\textsf{I}^2 },
\eequ
which describes the amount of $\Ldue(\mathcal{X},G_\bX)$ energy which is increased in the projection of $w(\cdot)$ onto $\mathcal{S}_\Lambda$ when the element \textsf{I} is refined. It also accounts for the decreased projection error when \textsf{I} is refined. If we were in a classical situation of Lebesgue measure and dyadic refinement, then $\nu^2_\textsf{I}(w)$ would
be exactly the sum of squares of the (scaling) Haar coefficients of $w(\cdot)$ corresponding to \textsf{I}.

We can use $\nu_\textsf{I} (w)$ to generate an adaptive partition. Given any $\lambda > 0$, let $\mathcal{T}(w,\lambda)$ be the smallest proper tree that contains all $\textsf{I} \in \mathcal{T}$ for which $\nu_\textsf{I} (w) \geqslant \lambda$. This tree can also be described as the set of all $\textsf{J} \in \mathcal{T}$ such that there exists $\textsf{I} \subset \textsf{J}$ which verifies $\nu_\textsf{I} (w) \geqslant \lambda$. Note that since $w \in \Ldue(\mathcal{X},G_\bX)$, the set of nodes such that $\nu_\textsf{I} (w) \geqslant \lambda$ is always finite and so is $\mathcal{T}( w ,\lambda)$. Corresponding to this tree we have the partition $\Lambda(w,\lambda)$ consisting of the outer leaves of $\mathcal{T}(w,\lambda)$. We shall define some new approximation spaces $\Beso^s$ which measure the regularity of a given function $w(\cdot)$ by the size of the tree $\mathcal{T}(w,\lambda)$.

Given $s > 0$, we let $\Beso^s$ be the collection of all $w \in \Ldue(\mathcal{X},G_\bX)$ such that the following is finite
\bequ\label{eq:CH3:binev_(5)}
|w|_{\Beso^s }^p  = \mathop {\sup }\limits_{\lambda  \geqslant 0} \big\{ {\lambda ^p \, \card\left( {\mathcal{T}(w,\lambda )} \right)} \big\},\quad \text{with}\;p = (s + \tfrac{1}
{2})^{ - 1}.
\eequ
We obtain the norm for $\Beso^s$ by adding $\| w \|_{\Ldue(G_\bX)}$ to $| w |_{\Beso^s}$. One can show that
\bequ\label{eq:CH3:binev_(6)}
\|w - \Pi _{\Lambda (w,\lambda )} (w)\|_{\Ldue(G_{\bX} )}  \leqslant C(s)\,|w|_{\Beso^s}^{\tfrac{1}{{2s + 1}}} \,\lambda ^{\tfrac{{2s}}{{2s + 1}}} \leqslant
    C(s)\,|w|_{\Beso^s } N^{ - s} ,
\eequ
where $N = \card(\mathcal{T}(w,\lambda ))$ and the constant $C(s)$ depends only on $s$ (see Cohen \emph{et al.} \cite{Cohen:Dah:Daub:deVore:2001}). It follows that every function $w \in \Beso^s$ can be approximated to order $\mathcal{O}(N^{-s})$ by $\Pi_\Lambda(w)(\cdot)$ for some partition $\Lambda$ with $\card(\Lambda) = N$. This should be contrasted with $\mathcal{A}^s$ which has the same approximation order for the uniform partition.
It is easy to see that $\Beso^s$ is larger than $\mathcal{A}^s$. In classical settings, the class $\Beso^s$ is well understood. For example, in the case of Lebesgue measure and dyadic partitions we know that each Besov space $\Beso^{\tau,s}_q$ with $\tau > (s/d + 1/2)^{-1}$ and $q\in(0,\infty]$ arbitrary, is contained in $\Beso^{s/d}$ (see \cite{Cohen:Dah:Daub:deVore:2001}). This should be compared with the $\mathcal{A}^s$ where we know that $\mathcal{A}^{s/d} = \Beso^{2,s}_\infty$ as we have noted earlier. In the next section we will see how to ``visualize'' these approximation spaces when we use \emph{warped} wavelet bases to build our partitions.

Until now, we have only considered the problem of approximating elements of some smoothness class by approximators associated to (adaptive) partitions of their domain $\mathcal{X}$: no data, no noise; just functions. Here, instead, we assume that $f(\cdot)$ denotes, as before, the regression function and we return to the problem of estimating it from a given data-set. Clearly, we can use the functions in $\mathcal{H} = \mathcal{S}_\Lambda$ for this purpose, so that the ``incarnation'' in this context of what we called the \emph{empirical minimizer}, is given by
\[
f_{\bz,\Lambda }  = \mathop {\arg \min }\limits_{w \in \mathcal{S}_\Lambda  } \frac{1}{n}\sum\limits_{i = 1}^n {[w(\bx_i ) - y_i ]^2 },
\]
the orthogonal projection of $y = y(\bx)$ onto $\mathcal{S}_\Lambda$ with respect to the empirical norm
\[
\| y \|_{\Ldue(\mathcal{X},\delta_\bX )}^2 =\frac{1}{n}\sum\limits_{i = 1}^n {|y(\bx_i )|^2 },
\]
with $y(\bx_i) = y_i$, and we can compute it by solving $\textsf{card}(\Lambda)$ independent problems, one for each element $\textsf{I} \in \Lambda$. The resulting estimator can than be written as
\[
f_{\bz,\Lambda } ( \cdot ) = \sum\limits_{\MyI \in \Lambda } {s_\MyI (\bz)
\tfrac{1}{{\sqrt {\mathbb{G}_{\bX,n} (\MyI)} }}} \boldsymbol{1}_\MyI ( \cdot ),
\]
where, for each $\MyI \in \Lambda$,
\[
s_\MyI (\bz) = \frac{1}{n}\sum\limits_{i = 1}^n
    {y_i \tfrac{1}{{\sqrt {\mathbb{G}_{\bX,n} (\MyI)}}}\boldsymbol{1}_\MyI (\bx_i )}
\quad \text{and} \quad
    \mathbb{G}_{\bX,n} (\MyI) = \frac{1}{n}\sum\limits_{i = 1}^n {\boldsymbol{1}_\MyI (\bx_i )},
\]
are the empirical counterparts of the theoretical coefficients defined in Equation (\ref{eq:CH3:piecewiseConsts_coeff}).
\begin{table}
\centering
\begin{tabular}{llllllllllll} \hline
{\textsf{{\textbf{Algorithm:} \textsf{Least--Squares on Adaptive Partitions}}}} \\ \hline
\textbf{Require:} Sample $\bz = \{(\bx_i,y_i)\}_{i\in\nSeq}$; threshold $\lambda_n$, $\gamma > 0$ smoothness index \\
\textbf{Output:} An estimator $f_\bz(\cdot)$ for the regression function $f(\cdot)$ \\
\textbf{Setup:} \\
$\;\;\; 1:$ Define $\MyJ^\star = \min\big\{j\in\MyNN:2^j\leqslant\lambda_n^{-1/\gamma}\big\}$ \\
\textbf{Generator:} \\
$\;\;\; 2:$ Compute $\nu_\MyI(\bz)$ for the nodes $\MyI$ at a refinement level $j<\MyJ^\star$\\
$\;\;\; 3:$ Threshold $\{\nu_\MyI(\bz)\}_{\MyI}$ at level $\lambda_n$ obtaining the set: $\Sigma(\bz,n) = \{\MyI \in \mathcal{T}_{\MyJ^\star}:\nu_\MyI(\bz) \geqslant \lambda_n\}$ \\
$\;\;\; 4:$ Complete $\Sigma(\bz,n)$ to a tree $\mathcal{T}(\bz,n)$ by adding nodes $\textsf{J} \supset \MyI \in \Sigma(\bz,n)$ \\
$\;\;\; 7:$ \textbf{Return} The estimator $f_\bz(\cdot)$ that minimizes the empirical risk on $\Lambda(\bz,n)$ \\
\hline
\end{tabular}
\caption[Algorithm: \textsf{Least--Squares on Adaptive Partitions}]{\textsf { \small Least--squares on adaptive partitions driven by the empirical residuals $\nu_\MyI(\bz)$ defined in Equation (\ref{eq:CH3:Empirical_residuals}). Adapted from \cite{Binev:Cohen:Dahmen:DeVore:Temlyakov:2004}.}}\label{algo:CH3:LS_on_AdaptParts}
\end{table}
With the coefficients $\{s_\MyI (\bz)\}_{\MyI\in\Lambda}$ at hand, we can build linear estimators $f_\bz(\cdot)$ corresponding to uniform partitions with cardinality suitably chosen to balance the bias and variance of $f_\bz(\cdot)$ when the true regression function $f(\cdot)$ belongs to some specific smoothness class. Alternatively, defining the empirical versions of the residuals introduced in Equation (\ref{eq:CH3:Theoretical_residuals}) as
\bequ\label{eq:CH3:Empirical_residuals}
\nu _{\MyI}(\bz) = \sqrt {\sum\limits_{\textsf{J} \in \mathcal{C}(\MyI)}
    {s^2_\textsf{J}(\bz)}  - s^2_\MyI (\bz)},
\eequ
we can mimic the adaptive procedure introduced in the previous section (see Table \ref{algo:CH3:LS_on_AdaptParts}) to get universal\footnote{A synonymous of ``adaptive'': the estimator does not require any prior knowledge of the smoothness of the regression function $f(\cdot)$.} estimators based on adaptive partitions. These partitions have the same tree structure as those used in the \textsf{CART} algorithm \cite{Breiman:Freidman:Olshen:Stone:1984}, yet the selection or the right partition is quite different since it is not based on an optimization problem but on a thresholding technique applied to to empirical quantities computed at each node of the tree which play a role similar to wavelet coefficients as we will see in the following (see \cite{Dono:1997} for a connection between \textsf{CART} and thresholding in one or several orthonormal bases).

\subsection{A Universal Algorithm Based on Warped Wavelets}\label{subsubsec:CH3:UniversalWarpedWavelet}

The choice we made in the previous Section of adopting piecewise constant functions as approximators, severely limits the optimal convergence rate to approximation spaces corresponding to smoothness classes of low or no pointwise regularity (see \cite{Binev:Cohen:Dahmen:DeVore:2007} for an interesting extension based on piecewise polynomial approximations). A possible way to fix this problem would be to use the complexity regularization approach for which optimal convergence results could be obtained in the piecewise polynomial context (see for instance Theorem 12.1 in \cite{Gyorfi:etal:2002}, and the paper by Kohler \cite{Kohler:2003}).

In the present context where the marginal design distribution $G_X(\cdot)$ is assumed to be known, we have another option based on the \emph{warped} systems introduced in Section \ref{sec:CH3:WarpedWavelet}.

It is worth mentioning that in this section we will concentrate on the $\mathcal{X} \equiv [0,1]$. The present setting could be generalized to the case where $G_X(\cdot)$ is a $d$--dimensional tensor product. However, the full generalization to dimension $d > 1$ is more involved and will not be discussed here.

To ``translate'' the concepts highlighted in the previous two sections in terms of \emph{warped} systems, consider a compactly supported wavelet basis $\{\psi_{j,k}(\cdot), j \geqslant -1, k\in\MyZ\}$, where $\psi_{-1,k}(\cdot) = \phi_{0,k}(\cdot)$ denotes the scaling function, and its \emph{warped} version $\{\psi_{j,k}(G_X(\cdot)), j \geqslant -1, k\in\MyZ\}$. Then, for each $f \in \Ldue([0,1],G_X)$, consider its expansion in this basis
\[
f(x) = \sum\limits_{j,k} {d_{j,k} \psi _{j,k} (G_X (x))}.
\]
In this context, a \emph{tree} is a finite set $\mathcal{T}$ of indexes $(j,k)$, $j \in \MyNN_0$ and $k \in \{0,\ldots,2^j-1\}$, such that $(j,k) \in \mathcal{T}$ implies $(j-1, \lfloor k/2 \rfloor) \in \mathcal{T}$, i.e., all ``ancestors'' of the point $(j,k)$ in the dyadic grid also belong to the tree.

One can then consider the best \emph{tree--structured} approximation to $f(\cdot)$, by trying to minimize
\[
\left\| {f - \sum\limits_{(j,k) \in \mathcal{T}} {d_{j,k} \psi _{j,k} (G)} } \right\|_{\Ldue(G_X )}^2,
\]
over all tree $\mathcal{T}$ having the same cardinality $N$, and all choices of $d_{j,k}$. However the procedure of selecting the optimal tree is costly in computational time, in comparison to the simple reordering that characterize the classical thresholding procedure described in the previous section. A more reasonable approach is to use suboptimal tree selection algorithms inspired by the adaptive procedure introduced before. In detail, we start from an initial tree $\mathcal{T}_0 = \{(0,0)\}$ and let it ``grow'' as follow:
\begin{enumerate}
\item Given a tree $\mathcal{T}_N$, define its ``leaves'' $\mathcal{L}(\mathcal{T}_N)$ as the indexes $(j,k) \notin \mathcal{T}_N$ such that $(j-1,\lfloor k/2 \rfloor) \in \mathcal{T}_N$.

\item For $(j,k) \in \mathcal{L}(\mathcal{T}_N)$ define the residual
\[
\nu_{j,k}  = \sqrt{ {\sum\limits_{\textsf{I}_{\ell ,m}  \subset \textsf{I}_{j,k} }
    {|d_{\ell ,m} |^2 } } },
\]
with $\textsf{I}_{j,k} = [2^{-j}k,2^{-j}(k+1)]$.

\item Choose $(j_0,k_0) \in \mathcal{L}(\mathcal{T}_N)$ such that
\[
\nu_{j_0,k_0} = \mathop {\max }\limits_{(j,k) \in \mathcal{L}(\mathcal{T}_N )} \nu_{j,k},
\]

\item Define $\mathcal{T}_{N+1} = \mathcal{T}_{N} \cup \{(j_0,k_0)\}$.
\end{enumerate}
Note that this algorithm can either be controlled by the cardinality $N$ of the tree, or by the size of the residuals as in Table \ref{algo:CH3:LS_on_AdaptParts}.

Now, let $\Lambda$ be the dyadic partition associated to any such tree, and define
\[
\Pi _\Lambda  (f)(x) = \sum\limits_{\textsf{I} \in \Lambda } {d_\textsf{I}
    \psi _\textsf{I} (G(x))} ,\quad \text{with} \quad d_\textsf{I}  =
    \left\langle {f,\psi _\textsf{I} (G)} \right\rangle _{\Ldue(G_X )},
\]
and their empirical counterparts
\[
f_{\Lambda ,\bz} (x) = \sum\limits_{\textsf{I} \in \Lambda } {d_\textsf{I} (\bz)
    \psi _\textsf{I} (G(x)),\quad \text{with} \quad d_\textsf{I} (\bz) =
    \frac{1}{n}\sum\limits_{i = 1}^n {Y_i \psi _\textsf{I} (G(X_i ))} }.
\]
Then, by adapting the techniques used in \cite{Binev:Cohen:Dahmen:DeVore:Temlyakov:2004}, in Section \ref{subsec:CH3:ProofUniform} we prove the following result for uniform partitions:
\begin{teo}\emph{(Optimality for uniform partitions)}\label{teo:CH3:uniformParts}
Assume that $f\in\mathcal{A}^s$ and define the estimator $f_\bz = f_{\bz,\Lambda_{\MyJ^\star}}$, with
\[
\MyJ^\star = \MyJ^\star(n) = \min \big\{ j \in \MyNN :  2^{j\,(1+2s)} \geqslant \tfrac{n}{\log(n)} \big\}.
\]
Then, given any $\beta > 0$, there is a constant $\tilde{c}$ such that
\[
\Prob^{ \otimes n} \left\{ {\left\| {f - f_\bz } \right\|_{\Ldue(G_X )}^2 > (\tilde c
    + |f|_{\mathcal{A}^s } )\left( {\tfrac{{\log n}}{n}} \right)^{\tfrac{s}{{2s + 1}}} }
    \right\} \leqslant C n^{ - \beta }
\]
and
\[
\Exp^{\otimes n} \left\{ {\left\| {f - f_\bz } \right\|_{\Ldue(G_X )}^2 } \right\} \leqslant
    (C + |f|_{\mathcal{A}^s }^2 )\left( {\tfrac{{\log n}}{n}} \right)^{\tfrac{{2s}}{{2s + 1}}},
\]
where $C$ depends only on $M$.
\end{teo}

Theorem \ref{teo:CH3:uniformParts} is satisfactory in the sense that the rate $\left[ {\tfrac{{\log (n)}}{n}} \right]^{ - s/(2s + 1)}$ is known to be optimal (or minimax) over the class $\mathcal{A}^s$ save for the logarithmic factor. However, it is unsatisfactory in the sense that the estimation procedure requires a--priori knowledge of the smoothness parameter $s$ which appears in the choice of the resolution level $j$. Moreover, as noted before, the smoothness assumption $f \in \mathcal{A}^s$ is too severe. Consequently, our next task, will consist in deriving a method capable of treating both defects. To this end, mimicking Equation (\ref{eq:CH3:Empirical_residuals}), we define the empirical residuals as
\[
\nu_{j,k} (\bz) = \sqrt{ {\sum\limits_{\textsf{I}_{\ell ,m}  \subset \textsf{I}_{j,k} }
    {|d_{\ell ,m} (\bz)|^2 } }  }.
\]
Then, for some\footnote{$\kappa$ is essentially a smoothing parameter to be selected by cross--validation, for instance). Notice that in our theoretical developments we will only assume that $\kappa$ is ``large enough'' to ensure the desired concentration inequalities.} $\kappa > 0$, let
\[
\lambda_n = \kappa \sqrt{\frac{\log(n)}{n}},
\]
be a given threshold. Now, adapting the algorithm given in Table \ref{algo:CH3:LS_on_AdaptParts}, assume that the estimator $f_\bz(\cdot)$ is generated as detailed in Table \ref{algo:CH3:warpedTree}.
\begin{table}
\centering
\begin{tabular}{llllllllllll} \hline
{\textsf{{\textbf{Algorithm:} \textsf{``Treed'' Approximations from \emph{Warped} Wavelets}}}} \\ \hline
\textbf{Require:} Sample $\bz = \{(\bx_i,y_i)\}_{i\in\nSeq}$; threshold $\lambda_n$, $\gamma \geqslant \tfrac{1}{2}$ smoothness index \\
\textbf{Output:} An estimator $f_\bz(\cdot)$ for the regression function $f(\cdot)$ \\
\textbf{Setup:} \\
$\;\;\; 1:$ Define $\MyJ^\star=\min\big\{j\in\MyNN:2^j\leqslant\lambda_n^{-1/\gamma}\big\}$ \\
\textbf{Generator:} \\
$\;\;\; 2:$ Compute $\nu_{j,k}(\bz)$ for any node $(j,k)$ at a refinement level $j<\MyJ^\star$\\
$\;\;\; 3:$ Threshold $\{\nu_{j,k}(\bz)\}_{j,k}$ at level $\lambda_n$ obtaining: $\Sigma(\bz,n) = \{(j,k) \in \mathcal{T}_{\MyJ^\star}:\nu_{j,k}(\bz) \geqslant \lambda_n\}$ \\
$\;\;\; 4:$ Complete $\Sigma(\bz,n)$ to a tree $\mathcal{T}(\bz,n)$ by adding nodes $(\ell,m)\in\mathcal{P}\big(\{(j,k)\}\big)$\\
$\;\;\;\;\;$ for all $(j,k)\in \Sigma(\bz,n)$ \\
$\;\;\; 7:$ \textbf{Return} The estimator $f_{\bz} ( \cdot ) = \sum\nolimits_{(j,k) \in \Lambda (\bz,n)} {d_{j,k} (z)\psi _{j,k} (G_X ( \cdot ))}$ \\
\hline
\end{tabular}
\caption[Algorithm: \textsf{``Treed'' Approximations from \emph{Warped} Wavelets}]{\textsf { \small Tree--structured approximations from \emph{warped} wavelet decompositions.}}\label{algo:CH3:warpedTree}
\end{table}
\noindent Then, in Section \ref{subsec:CH3:ProofAdaptive}, we prove the following

\begin{teo}\emph{(Optimality for ``growing'' adaptive partitions)}\label{teo:CH3:adaptiveParts}
Let $\beta$ and $\gamma \geqslant \tfrac{1}{2}$ be arbitrary. Then, there exists $\kappa > 0$, such that, whenever $f \in \mathcal{A}^\gamma \cap \Beso^s$ for some $s>0$, the following inequalities hold
\[
\Prob^{ \otimes n} \left\{ {\left\| {f - f_\bz } \right\|_{\Ldue(G_X )}  > \tilde c\left( {\tfrac{{\log n}}{n}} \right)^{\tfrac{s}{{2s + 1}}} } \right\} \leqslant Cn^{ - \beta },
\]
and
\[
\Exp^{ \otimes n} \left\{ {\left\| {f - f_\bz } \right\|_{\Ldue(G_X )}^2 } \right\} \leqslant C\left( {\tfrac{{\log n}}{n}} \right)^{\tfrac{{2s}}{{2s + 1}}},
\]
where the constants $\tilde{c}$ and $C$ do not depend on the sample size $n$.
\end{teo}
Theorem \ref{teo:CH3:adaptiveParts} is definitively more satisfactory than Theorem \ref{teo:CH3:uniformParts} in two respects:
\begin{itemize}
\item The optimal rate $\left[ {\tfrac{{\log (n)}}{n}} \right]^{ - s/(2s + 1)}$ is now obtained under weaker smoothness assumptions on the regression function, namely, $f \in \Beso^s$ in place of $f \in \mathcal{A}^s$, with the extra assumption $f \in \mathcal{A}^\gamma$ with $\gamma \geqslant \tfrac{1}{2}$ arbitrary.
\item The estimator we obtain is adaptive (universal), in the sense that the value of $s$ does not enter the definition of the algorithm. The procedure automatically extract information about the regularity of the regression function from the data at hand.
\end{itemize}
It is interesting to notice that in standard thresholding (standard denoising or density estimation, for instance) one usually sets the highest level $\MyJ^\star$ so that $2^{\MyJ^\star} \sim n/\log(n)$; here we have to stop much sooner, namely, $2^{\MyJ^\star} \sim \sqrt{n/\log(n)}$, as in \cite{KP04}. This is especially necessary to obtain the exponential inequalities in Section \ref{subsec:CH3:ProofUniform} and \ref{subsec:CH3:ProofAdaptive}.

A final remark on the approximation spaces $\mathcal{A}^s$ and $\Beso^s$ is in order. In a previous section, we mentioned that, when $G_X(\cdot)$ is the Lebesgue measure, then the spaces $\mathcal{A}^s$ and $\Beso^s$ are well understood. In particular, each Besov space $\Beso^{\tau,s}_q$ with $\tau > (s + 1/2)^{-1}$ and $q \in (0,+\infty]$, is contained in $\Beso^s$ (see Cohen \emph{et al.} \cite{Cohen:Dah:Daub:deVore:2001,Cohen:Daub:Gul:Orch:2002}), whereas $\mathcal{A}^s = \Beso^{2,s}_\infty$. For general partitions it is not totally clear how to express the content of these approximation spaces in terms of reasonably simple variations of common smoothness classes. Things get slightly simpler when we employ a \emph{warped} wavelet basis to generate the partition. As a result, we can map approximation properties imposed on $f(\cdot)$ to regularity properties over its warped version $f\circ G_X^{-1}(\cdot)$. So, assuming $f \in \mathcal{A}^s$ is equivalent to impose $f\circ G^{-1}_X \in \Beso^{2,s}_\infty$ as soon as $\mathcal{A}^s$ is defined in terms of warped wavelets.

\section{Discussion}\label{sec:CH3:Discu}

The dependence on the design marginal $G_X(\cdot)$ is so far a clear weakness of our approach from both a theoretical and a practical perspective. Nevertheless, an obvious option to extend our tree--structured procedure to the case of an unknown $G_X(\cdot)$, would probably end up combining the arguments introduced in \cite{Binev:Cohen:Dahmen:DeVore:Temlyakov:2004} with those considered by Kerkyacharian and Picard in \cite{KP04} and \cite{Kerk:Picard:2005}. Another practical option might be to adopt a \emph{split sample} approach and measure smoothness in terms of the \emph{discrete} norm induced by the data. Here we also mention the fact that in Theorem \ref{teo:CH3:adaptiveParts} we require the knowledge of the parameter $\gamma$ which can be arbitrary close to $1/2$. As in \cite{Binev:Cohen:Dahmen:DeVore:Temlyakov:2004}, it is probably possible to remove the dependency on $\gamma$ at the price of using the much more complicated construction proposed by Binev and DeVore in \cite{Binev:deVore:2004}.

\section{Proofs for Section \ref{subsec:CH3:TreeApproximation}}\label{sec:CH3:ProofsWarped}

\subsection{Proof of Theorem \ref{teo:CH3:uniformParts}}\label{subsec:CH3:ProofUniform}

For any given partition $\Lambda$, a natural way to control $\|f - f_{\bz,\Lambda}\|^2_{\Ldue(G_X)}$ is by splitting it into a bias and variance term denoted respectively with $e_1$ and $e_2$ in the following equation
\bequ\label{eq:CH3:control_basic}
\|f - f_{\bz,\Lambda }\|_{\Ldue(G_X)}^2  = \|f - \Pi _\Lambda(f)\|_{\Ldue(G_X )}^2
    + \|\Pi _\Lambda (f) - f_{\bz,\Lambda } \|_{\Ldue(G_X )}^2  = e_1  + e_2 .
\eequ
$e_1$ will be controlled by using the smoothness assumptions we made in the statement of the theorem, whereas the variance term $e_2$ will be controlled by Bernstein's inequality.

Lets start with this second step observing that, by denoting $[d_\textsf{I}  - d_\textsf{I} (\bz)]$ with $\Delta_\textsf{I}(\bz)$, then by orthonormality of the warped system we have
\[
\|\Pi _\Lambda  (f) - f_{\bz,\Lambda } \|_{\Ldue(G_X )}^2  =
    \biggl\| \sum\nolimits_{\textsf{I}\in\Lambda} {[d_\textsf{I}  - d_\textsf{I} (\bz)]} \,
    \psi_\textsf{I}(G_X(\cdot))\biggr\|_{\Ldue(G_X )}^2  =
    \sum\nolimits_{\textsf{I}\in\Lambda} {\Delta _\textsf{I}^2 (\bz)}.
\]
Hence, for any $\eta > 0$,
\begin{eqnarray}\label{eq:CH3:control_var_1}
& & \Prob^{\otimes n}\left\{ {\|\Pi _\Lambda(f) - f_{\bz,\Lambda}\|_{\Ldue(G_X)}>\eta} \right\}
    = \Prob^{ \otimes n} \left\{ {\sum\nolimits_{\textsf{I} \in \Lambda }
    {\Delta _\textsf{I}^2 (\bz)} > \eta ^2 } \right\} \leqslant \\
& & \quad\quad\quad\quad\quad
    \leqslant \card(\Lambda )\cdot\Prob^{ \otimes n} \left\{ {\Delta _\textsf{I}^2 (\bz) >
    \tfrac{{\eta ^2 }}{{\card(\Lambda )}}} \right\} = \nonumber \\
& & \quad\quad\quad\quad\quad
    =  \card(\Lambda) \cdot \Prob^{ \otimes n} \biggl\{ {|\Delta _\textsf{I} (\bz)| >
    \tfrac{\eta }{{\sqrt {\card(\Lambda )} }}} \biggr\}, \nonumber
\end{eqnarray}
Consequently to control $e_2$ we just need to control $|\Delta_{\textsf{I}}(\bz)|$ and the cardinality of $\Lambda$. Now, if we define $U = Y \psi_{j,k}(X)$, then
\begin{itemize}
\item $M = \|U - \Exp\{U\}\|_\infty   \leqslant 2 \cdot 2^{j/2} \|\psi\|_\infty  \|f\|_\infty,$
\item $\sigma ^2  = \Exp\big\{ \, |U - \Exp(U)|^2 \big\}  \leqslant \Exp\big\{ \, |U|^2 \big\} \leqslant \|f\|_\infty ^2,$
\end{itemize}
as
\begin{eqnarray*}
\Exp\left\{ {\,|\psi _{j,k}(G(X))|^2 } \right\} \!\!\!\!&=& \!\!\!\!
    \int {|\psi _{j,k} (G(x))|^2 \mathrm{d}G_X (x)}  =
    \int\limits_0^1 {\big|\psi _{j,k} \big(G(G^{ - 1} (y))\big)\big|^2 \mathrm{d}y} = \\
    \!\!\!\! &=& \!\!\!\! \int\limits_0^1 {|\psi _{j,k} (y)|^2 \mathrm{d}y}  = 1.
\end{eqnarray*}
Hence, for any $\eta > 0$, by Bernstein's inequality we get

\begin{footnotesize}
\begin{eqnarray}\label{eq:CH3:control_var_2}
\Prob^{ \otimes n} \left\{ {|\Delta _\textsf{I} (\bz)| \geqslant
    \tfrac{\eta }{{\sqrt {\card(\Lambda )} }}} \right\} &\leqslant&
    2\exp \left\{ { - \frac{{n\,\eta ^2 }}{{2\,\card(\Lambda )
    \left[ {\|f\|_\infty ^2  + \tfrac{{2 \|\psi \|_\infty  \|f\|_\infty  }}{3}\tfrac{{2^{j/2} }}
    {{\sqrt {\card(\Lambda )} }}\eta } \right]}}} \right\} \leqslant \nonumber \hfill \\
&\leqslant&
    2\exp \left\{ { - \frac{{3n\,\eta ^2 }}{{C' \card(\Lambda )(3 + \eta )}}} \right\},
\end{eqnarray}
\end{footnotesize}

\noindent where $C' = 2 \max\{\|f\|^2_\infty,2\|\psi\|_\infty \|f\|_\infty\}$, and the last inequality comes from the fact that for any $\textsf{I} \in \Lambda$ we have $2^j  = |\textsf{I}|^{ - 1} \leqslant \card(\Lambda ) = 2^\MyJ$ for some $\MyJ \in \MyNN$, being $\Lambda$ a dyadic partition.

Now, back to our specific case. First of all remember that, by definition,
\[
\MyJ^\star = \MyJ^\star(n) = \min \big\{ j \in \MyNN :  2^{j\,(1+2s)} \geqslant \tfrac{n}{\log(n)} \big\},
\]
so
\bequ\label{eq:CH3:control_bias_2}
\card(\Lambda _{\MyJ^{\star} } ) \leqslant 2^{\MyJ^{\star}  + 1}  \leqslant
    2^2 2^{\MyJ^{\star}  - 1} \leqslant 2^2 \left[ {\tfrac{{\log(n)}}{n}}
    \right]^{-\tfrac{1}{{1 + 2s}}}.
\eequ
Hence, by definition of $\mathcal{A}^s$, we get the following bound for $e_1$:
\bequ\label{eq:CH3:control_bias}
\|f - \Pi _{\Lambda _{\MyJ^{\star} } } (f)\|_{\Ldue(G_X )}  \leqslant |f|_{\mathcal{A}^s }
    2^{ - \MyJ^{\star} s} \leqslant |f|_{\mathcal{A}^s } \left[ {\tfrac{{\log (n)}}{n}}
    \right]^{\tfrac{s}{{1 + 2s}}}.
\eequ
From Equation (\ref{eq:CH3:control_basic}) we then get
\[
\|f - f_{\bz,\Lambda _{\MyJ^{\star} } } \|_{\Ldue(G_X )}^2  \leqslant |f|_{\mathcal{A}^s }^2
    \left[ {\tfrac{{\log (n)}}{n}} \right]^{\tfrac{{2s}}{{1 + 2s}}}  + e_2^2 ,
\]
therefore, for all $\delta > 0$
\[
\Prob^{ \otimes n} \left\{ {\|f - f_{\bz,\Lambda _{\MyJ^{\star} } } \|_{\Ldue(G_X )}
    \geqslant \delta } \right\} \leqslant
    \Prob^{ \otimes n} \left\{ {e_2  > \delta  - |f|_{\mathcal{A}^s }
    \left[ {\tfrac{{\log (n)}}{n}} \right]^{\tfrac{{s}}{{1 + 2s}}} } \right\}.
\]
Taking $\delta = (\tilde c + |f|_{\mathcal{A}^s } )\left[ {\tfrac{{\log (n)}}{n}} \right]^{\tfrac{s}{{2s + 1}}}$ as in the statement of Theorem \ref{teo:CH3:uniformParts}, and applying Equations (\ref{eq:CH3:control_var_1}) and (\ref{eq:CH3:control_var_2}) noticing that $\left[ {\tfrac{{\log (n)}}{n}} \right]^{\tfrac{s}{{2s + 1}}} < 1$ for every $s>0$, we obtain

\begin{footnotesize}
\[
\Prob^{ \otimes n} \left\{ {e_2  > \tilde{c} \left[ {\tfrac{{\log (n)}}{n}}
    \right]^{\tfrac{s}{{1 + 2s}}} } \right\}
    \leqslant 2 \card(\Lambda_{\MyJ^{\star} })
    \exp \left\{ { - \left( {\tfrac{n}{{\card(\Lambda_{\MyJ^{\star} } )}}
    \left[ {\tfrac{{\log (n)}}{n}} \right]^{\tfrac{{2s}}{{1 + 2s}}} } \right)
    \tfrac{{3\,\tilde c^2 }}{{C'  (3 + \tilde{c})}}} \right\}.
\]
\end{footnotesize}

\noindent But from Equation (\ref{eq:CH3:control_bias_2}) we know how to bound the cardinality of our partition, therefore

\begin{footnotesize}
\be
\left( {\tfrac{n}{{\card(\Lambda _{\MyJ^{\star} } )}}\left[ {\tfrac{{\log (n)}}{n}} \right]^{\tfrac{{2s}}{{1 + 2s}}} } \right)\frac{{3\,\tilde c^2 }}
{{C'  (3 + \tilde c)}} &>& \left( {\tfrac{n}{{2^2 }}
\left[ {\tfrac{{\log (n)}}{n}} \right]^{\tfrac{1}{{1 + 2s}}}
\left[ {\tfrac{{\log (n)}}{n}} \right]^{\tfrac{{2s}}{{1 + 2s}}} } \right)
\frac{{3\,\tilde c^2 }}{{C'  (3 + \tilde c)}} = \\
&=& g(\tilde c) \cdot \log (n),
\ee
\end{footnotesize}

\noindent with
\[
g(\tilde c) = \frac{{3\,\tilde c^2 }}{{4 C'  (3 + \tilde c)}},
\]
so that
\be
\Prob^{ \otimes n} \left\{ {\|f - f_{\bz,\Lambda _{\MyJ^{\star} } } \|_{\Ldue(G_X )}
    \geqslant \delta } \right\} &=& \Prob^{ \otimes n} \left\{ {e_2  >
    \tilde c\left[ {\tfrac{{\log (n)}}{n}} \right]^{\tfrac{s}{{1 + 2s}}} } \right\} \leqslant \\
   &\leqslant& 2 \cdot 2^2 \left[ {\tfrac{{\log (n)}}{n}} \right]^{ - \tfrac{1}{{1 + 2s}}}
    \exp \left\{ {\log \left[ {n^{ - g(\tilde c)} } \right]} \right\} \leqslant  \hfill \\
   &\leqslant& 8\,n^{ - [g(\tilde c) - 1]}  \leqslant 8\,n^{ - \beta },
\ee
where the last inequality holds as soon as $g(\tilde c) - 1 > \beta$. And this complete the proof since from here we can easily derive a bound for the risk by using Equation (\ref{eq:CH3:Exp2Prob}).

\subsection{Proof of Theorem \ref{teo:CH3:adaptiveParts}}\label{subsec:CH3:ProofAdaptive}

Lets start with a bit of notation. First of all, for each $\lambda > 0$, we will denote by
\begin{itemize}
\item $\mathcal{T}(f,\lambda)$: smallest tree which contains all dyadic intervals $\textsf{I}$ such that $\nu_\textsf{I} > \lambda$.
\item $\Lambda(f,\lambda)$: partition induced by the outer leaves of $\mathcal{T}(f,\lambda)$.
\item $\mathcal{T}(f,\lambda,\bz)$: smallest tree which contains all dyadic intervals $\textsf{I}$ such that $\nu_\textsf{I}(\bz) > \lambda$.
\item $\Lambda(f,\lambda,\bz)$: partition induced by the outer leaves of $\mathcal{T}(f,\lambda,\bz)$.
\end{itemize}
If $\Lambda_0$ and $\Lambda_1$ are partitions associated to the tree $\mathcal{T}_0$ and $\mathcal{T}_1$, then we denote by
\begin{itemize}
\item $\Lambda_0 \vee \Lambda_1$ the partition associated to the tree $\mathcal{T}_0 \cup \mathcal{T}_1$,
\item $\Lambda_0 \wedge \Lambda_1$ the partition associated to the tree $\mathcal{T}_0 \cap \mathcal{T}_1$.
\end{itemize}
Finally, let $\lambda_n = \kappa \sqrt{\tfrac{\log(n)}{n}}$ for some $\kappa > 0$, and
\[
\MyJ^\star = \min\big\{j\in\MyNN : 2^j \leqslant \lambda_n^{-1/\gamma}\big\}.
\]
Then for each $\lambda > 0$, define the partitions
\begin{itemize}
\item $\Lambda(\lambda) = \Lambda(f,\lambda)\wedge\Lambda_{\MyJ^\star}$,
\item $\Lambda(\lambda,\bz) = \Lambda(f,\lambda,\bz)\wedge\Lambda_{\MyJ^\star}$.
\end{itemize}
Therefore, in this section, we consider the adaptive estimator
\[
f_{\bz,n}(x)  = f_{\bz,\Lambda (\lambda_n ,\bz)}(x) = \sum\limits_{\textsf{I}
    \in \Lambda (\tau _n ,\bz)} {d_\textsf{I} (\bz)\psi _\textsf{I} (G_X(x))}.
\]

Lets now start the proof observing that, using the triangle inequality, we can decompose the loss as follow
\[
\|f - f_{\bz,n} \|_{\Ldue(G_X )}  = e_1  + e_2  + e_3  + e_4 ,
\]
where
\begin{itemize}
\item $e_1  = \|f - \Pi _{\Lambda (\lambda _n ,\bz) \vee \Lambda (2\lambda _n )} (f)\|_{\Ldue(G_X )},$
\item $e_2  = \|\Pi _{\Lambda (\lambda _n ,\bz) \vee \Lambda (2\lambda _n )} (f)
    - \Pi _{\Lambda (\lambda _n ,\bz) \wedge \Lambda(2^{-1} \lambda _n )}(f)\|_{\Ldue(G_X )},$
\item $e_3  = \|\Pi _{\Lambda (\lambda _n ,\bz) \wedge \Lambda (2^{ - 1} \lambda _n )} (f) - f_{\bz,\Lambda (\lambda _n ,\bz) \wedge \Lambda (2^{ - 1} \lambda _n )} \|_{\Ldue(G_X )},$
\item $e_4  = \|f_{\bz,\Lambda (\lambda _n ,\bz) \wedge \Lambda (2^{ - 1} \lambda _n )}  - f_{\bz,\Lambda (\tau _n ,\bz)} |_{\Ldue(G_X )}.$
\end{itemize}
This type of splitting is frequently used in the analysis of wavelet thresholding procedures to deal with the fact that the partition built from those $\textsf{I}$ such that $\nu_{\textsf{I}}(\bz) \geqslant \lambda_n$, does not exactly coincides with the partition which would be chosen by an oracle based on those $\textsf{I}$ such that $\nu_\textsf{I} \geqslant \lambda_n$. This is accounted by the terms $e_2$ and $e_4$ which correspond to those dyadic interval $\textsf{I}$ such that $\nu_\textsf{I}(\bz)$ is significantly larger or smaller than $\nu_\textsf{I}$ respectively, and which will proved to be small in probability. The remaining terms $e_1$ and $e_3$ correspond respectively to the bias and variance of oracle estimators based on partitions obtained by zero--thresholding based on the unknown quantities $\{\nu_\textsf{I}\}_\textsf{I}$.

The first term $e_1$, being a bias, is treated by a deterministic estimate as in \cite{Binev:Cohen:Dahmen:DeVore:Temlyakov:2004}. More specifically, since $\Lambda (\lambda _n ,\bz) \vee \Lambda (2\lambda _n )$ is a refinement of $\Lambda (2 \lambda _n) = \Lambda (f,2 \lambda _n) \wedge \Lambda_{\MyJ^\star}$, we have (almost surely):
\be
e_1  &\leqslant&
    \|f - \Pi _{\Lambda (2\lambda _n )} (f)\|_{\Ldue(G_X )} \leqslant \\
&\leqslant&
    \|f - \Pi _{\Lambda (f,2\lambda_n )}(f)\|_{\Ldue(G_X )}  +
    \|\Pi _{\Lambda (f,2\lambda _n )}(f)-\Pi_{\Lambda (2\lambda _n )}(f)\|_{\Ldue(G_X )} \leqslant  \hfill \\
&\leqslant&
    \|f - \Pi _{\Lambda (f,2\lambda _n )}(f)\|_{\Ldue(G_X )} +
    \|f - \Pi _{\Lambda_{\MyJ^{\star}}}(f)\|_{\Ldue(G_X )}  \leqslant  \hfill \\
&\leqslant&
    C(s)[2\lambda _n ]^{\tfrac{{2s}}{{2s + 1}}} |f|_{\Beso^s}  +
    2^{ - \gamma \, \MyJ^{\star} } |f|_{\mathcal{A}^\gamma}\leqslant  \hfill \\
&\leqslant&
    C(s)[2\lambda _n ]^{\tfrac{{2s}}{{2s + 1}}} |f|_{\Beso^s}  +
    2^{ - \gamma } \lambda _n |f|_{\mathcal{A}^\gamma}.
\ee
Therefore
\[
e_1  \leqslant C(s)\left\{ {(2\kappa )^{\tfrac{{2s}}{{2s + 1}}}
    + 2^\gamma  \kappa } \right\}\max \left\{ {\,|f|_{\mathcal{A}^\gamma } ,|f|_{\Beso^s } } \right\} \left[ {\tfrac{{\log (n)}}{n}} \right]^{\tfrac{s}{{2s + 1}}}  \!\! = c_1 \left[
    {\tfrac{{\log (n)}}{n}} \right]^{\tfrac{s}{{2s + 1}}},
\]
as soon as $f \in \Beso^s \cap \mathcal{A}^\gamma$, with $c_1 = C(s)\left\{ {(2\kappa )^{\tfrac{{2s}}{{2s + 1}}} + 2^\gamma  \kappa } \right\}\max \left\{ {\,|f|_{\mathcal{A}^\gamma } ,|f|_{\Beso^s } } \right\}$.

The third term $e_3$ is treated by the estimate provided by combining Equations (\ref{eq:CH3:control_var_1}) and (\ref{eq:CH3:control_var_2})
\bequ\label{eq:CH3:e_3_bound_1}
\Prob^{ \otimes n} \left\{ {e_3  > \eta } \right\} \leqslant 2 \card(\Lambda _3 )
    \exp \left\{ { - \frac{{3n\,\eta ^2 }}{{C'
    \card(\Lambda _3 )(3 + \eta )}}} \right\},
\eequ
where $\Lambda_3 = \Lambda (\lambda _n ,\bz) \wedge \Lambda (2^{ - 1} \lambda _n )$. So
\begin{eqnarray}\label{eq:CH3:e_3_bound_2}
\card(\Lambda _3 ) &\leqslant&
    \card\left( {\Lambda (2^{ - 1} \lambda _n )} \right) =
    \card\left( {\Lambda (f,2^{ - 1} \lambda _n )
        \wedge \Lambda _{\MyJ^{\star} } } \right) \leqslant \nonumber \hfill \\
&\leqslant&
    \card\left( {\Lambda (f,2^{ - 1} \lambda _n )} \right) \leqslant
    (2^{ - 1} \lambda _n )^{ - p} \, |f|_{\Beso^s }^p  =
    2^p  {\lambda _n^{ - \tfrac{2}{{1 + 2s}}} } |f|_{\Beso^s }^p  = \nonumber \hfill \\
&=& 2^p \kappa ^{ - \tfrac{2}{{1 + 2s}}} \, |f|_{\Beso^s }^p
    \left[ {\tfrac{{\log (n)}}{n}} \right]^{ - \tfrac{2}{{2(1 + 2s)}}} =
    c_3 \left[ {\tfrac{{\log (n)}}{n}} \right]^{ - \tfrac{1}{{1 + 2s}}},
\end{eqnarray}
where we have used the fact that $1/p = 1/2 + s$.

For the remaining two terms, $e_2$ and $e_4$ we will show that $\forall \;\beta  > 0$ we fix, there exists a constant $C' > 0$ such that:
\bequ\label{eq:CH3:e_2&e_4_goal}
    \Prob^{ \otimes n} \left\{ {e_2  > 0} \right\} +
    \Prob^{ \otimes n} \left\{ {e_4 > 0} \right\}
    \leqslant C' \, n^{ - \beta }.
\eequ
Before we prove this, lets show why it is sufficient. Let $0 < \delta = \tilde{c} \left[\tfrac{\log(n)}{n}\right]^{\tfrac{1}{1+2s}}$ as in the statement of Theorem \ref{teo:CH3:adaptiveParts}. Then we have

\begin{footnotesize}
\be
\Prob^{ \otimes n} \left\{ {\|f - f_{\bz,n} \|_{\Ldue(G_X )} \geqslant \delta } \right\}
\!\!\!\!\!\!\!\!\!\! &\leqslant& \!\!\!\!\!\!\!\!\!\!
    \Prob^{ \otimes n} \left\{ {e_1  + e_2  + e_3  + e_4  \geqslant \delta } \right\}
    \leqslant \hfill \\
\!\!\!\!\!\!\!\!\!\! &\leqslant& \!\!\!\!\!\!\!\!\!\!
    \Prob^{ \otimes n} \left\{ {e_2  + e_3  + e_4  \geqslant
    (\tilde c - c_1 )\left[ {\tfrac{{\log (n)}}{n}} \right]^{\tfrac{s}{{2s + 1}}}}
    \right\} \leqslant  \hfill \\
\!\!\!\!\!\!\!\!\!\! &\leqslant& \!\!\!\!\!\!\!\!\!\!
    \Prob^{ \otimes n} \left\{ {e_2  > 0} \right\} +
    \Prob^{ \otimes n} \left\{ {e_4  > 0} \right\} +
    \Prob^{ \otimes n} \left\{ {e_3  \geqslant \tilde \delta } \right\} \leqslant  \hfill \\
\!\!\!\!\! &\mathop  \leqslant \limits^{\text{by Eq.(\ref{eq:CH3:e_2&e_4_goal})}}& \!\!\!\!\!
    C'n^{ - \beta } + \Prob^{ \otimes n} \left\{ {e_3 \geqslant \tilde \delta } \right\},
\ee
\end{footnotesize}

\noindent where $\tilde \delta  = (\tilde c - c_1 )\left[ {\tfrac{{\log (n)}}{n}} \right]^{\tfrac{s}{{2s + 1}}}$. Repeating the steps used to derive the bound we needed in Section \ref{subsec:CH3:ProofUniform}, from Equations (\ref{eq:CH3:e_3_bound_1}) and (\ref{eq:CH3:e_3_bound_2}), we obtain
\be
& & \left( {\tfrac{n}{{\card(\Lambda _3 )}}
    \left[ {\tfrac{{\log (n)}}{n}} \right]^{\tfrac{{2s}}{{1 + 2s}}} } \right)
    \frac{{3\,(\tilde c - c_1 )^2 }}{{C'  [3 + (\tilde c - c_1 )]}} > \\
& & \quad\quad\quad\quad\quad
    > \left( {\tfrac{n}{{c_3 }}\left[ {\tfrac{{\log (n)}}{n}} \right]^{\tfrac{1}{{1 + 2s}}}
    \left[ \tfrac{\log (n)}{n} \right]^{\tfrac{{2s}}{{1 + 2s}}} } \right)
    \frac{{3\,(\tilde{c}-c_1)^2 }}{{C'  [3 + (\tilde{c}-c_1)]}} =  \hfill \\
& & \quad\quad\quad\quad\quad
    = g(\tilde c) \cdot \log (n),
\ee
where
\[
g(\tilde c) = \frac{{3\,(\tilde{c}-c_1)^2 }}{{c_3 C'  [3 + (\tilde{c}-c_1)]}}.
\]
Therefore
\[
\Prob^{ \otimes n} \left\{ {e_3  \geqslant \tilde \delta } \right\} \leqslant
    2 c_3 \left[ {\tfrac{{\log (n)}}{n}} \right]^{ - \tfrac{1}{{1 + 2s}}} n^{ - g(\tilde c)}
    \leqslant
    c'n^{ - [g(\tilde c) - 1]}  \leqslant c'n^{ - \beta } ,
\]
as soon as $g(\tilde{c}) - 1 \geqslant \beta$. And this would conclude the proof of Theorem \ref{teo:CH3:adaptiveParts}.

We need to prove Equation (\ref{eq:CH3:e_2&e_4_goal}). The main tool is the following lemma
\begin{lem}\label{lem:CH3:residuals}
For each $\textsf{I} \in \Lambda_{\MyJ^\star}$, we have
\begin{itemize}
\item $\Prob^{ \otimes n} \left( {\{ \nu_\textsf{I}(\bz) \leqslant \lambda _n \}
    \cap \{ \nu _\textsf{I} \geqslant 2\lambda _n \} } \right)
    \leqslant 4 \,n^{ - g(\kappa )}$,
\item $\Prob^{ \otimes n} \left( {\{ \nu_\textsf{I}(\bz) \geqslant \lambda _n \}
    \cap \{ \nu_\textsf{I} \leqslant 2^{ - 1} \lambda _n \} } \right)
    \leqslant 4 \,n^{ - g(\kappa )}$,
\end{itemize}
where
\[
g(\kappa ) = \frac{{3\kappa ^2 }}{{8 C'\biggl( {3 +
    \kappa ^{1 - \tfrac{1}{{2\gamma }}} } \biggr)}}.
\]
\end{lem}
Before we prove Lemma \ref{lem:CH3:residuals}, lets show why this is sufficient. Remember that
\[
e_2  = \|\Pi _{\Lambda (\lambda _n ,\bz) \vee \Lambda (2\lambda _n )} (f) -
    \Pi _{\Lambda (\lambda _n ,z) \wedge \Lambda (2^{ - 1} \lambda _n )} (f)\|_{\Ldue(G_X )}.
\]
Consequently
\begin{itemize}
\item $e_2 \equiv 0$ if $\mathcal{T}(\lambda _n,\bz) \cup \mathcal{T}(2\lambda _n ) \equiv \mathcal{T}(\lambda _n ,\bz) \cap \mathcal{T}(2^{ - 1} \lambda _n )$,
\item $e_2 > 0$ if
\begin{scriptsize}
\be
\mathcal{T}(\lambda_n,\bz) \cup \mathcal{T}(2\lambda_n ) \supset \mathcal{T}(\lambda_n,\bz)
    \cap \mathcal{T}(2^{-1}\lambda_n)
    \!\!\!\! &\Leftarrow& \!\!\!\!
\bfun
   \mathcal{T}(\lambda_n ,\bz) \not\subset \mathcal{T}(2^{-1}\lambda_n)  \\
   \text{or}  \\
   \mathcal{T}(2\lambda_n ) \not\subset \mathcal{T}(\lambda_n ,\bz)
\efun
\\
\!\!\!\! &\Leftarrow& \!\!\!\! \exists \;\textsf{I}\;\text{s.t.}\;
\bfun
   \{ \nu_\textsf{I}(\bz) \leqslant \lambda_n \} \cap
    \{ \nu_\textsf{I} \geqslant 2\lambda _n \}  \\
   \text{or}  \\
   \{ \nu_\textsf{I}(\bz) \geqslant \lambda_n \} \cap
    \{ \nu_\textsf{I}  \leqslant 2^{ - 1} \lambda _n \}
\efun .
\ee
\end{scriptsize}
\end{itemize}
Therefore
\begin{eqnarray}\label{eq:CH3:e_2_A}
\Prob^{ \otimes n} \left\{ {e_2  > 0} \right\} &\leqslant&
    \sum\limits_{\textsf{I} \in \Lambda_{\MyJ^{\star} } }
    {\Prob^{ \otimes n} \left( {\{ \nu_\textsf{I} (\bz) \leqslant \lambda _n \}
    \cap \{ \nu_\textsf{I} \geqslant 2\lambda _n \} } \right)}  +  \hfill \\
& & \quad\quad
    + \sum\limits_{\textsf{I} \in \Lambda_{\MyJ^{\star} } }
    {\Prob^{ \otimes n} \left( {\{ \nu_\textsf{I}(\bz) \geqslant \lambda_n \} \cap
    \{\nu_\textsf{I} \leqslant 2^{ - 1} \lambda _n \} } \right)}  = R_1 + R_2. \nonumber
\end{eqnarray}
Then, by applying the first part of Lemma \ref{lem:CH3:residuals}, we get
\begin{eqnarray}\label{eq:CH3:e_2_B}
R_1  &\leqslant&
    \card(\Lambda_{\MyJ^{\star} } ) \, 4 \, n^{ - g(\kappa )}  \leqslant
    \card(\Lambda_0)\,2^{\MyJ^{\star}} 4 \, n^{ - g(\kappa )}  \leqslant \hfill \\
&\leqslant&
    \card(\Lambda_0 )\,\lambda_n^{-1/\gamma } 4 \,n^{-g(\kappa)} \leqslant
    \card(\Lambda _0 )\,\kappa ^{ - 1/\gamma } \left[ {\tfrac{n}{{\log (n)}}}
        \right]^{\tfrac{1}{{2\gamma }}} 4 \,n^{ - g(\kappa )}  \leqslant \hfill \nonumber \\
&\leqslant&
    \card(\Lambda _0 )\,\kappa ^{ - 1/\gamma } n^{\tfrac{1}{\gamma }}
    4 \, n^{ - g(\kappa )}  = C' \, n^{ -\left[g(\kappa ) - \tfrac{1}{\gamma }\right]}, \nonumber
\end{eqnarray}
and analogously, by the second part of Lemma \ref{lem:CH3:residuals}, we obtain
\bequ\label{eq:CH3:e_2_C}
R_2  \leqslant C'\,n^{ - \left[g(\kappa ) - \tfrac{1}{\gamma }\right]}.
\eequ

Applying again the second part of Lemma \ref{lem:CH3:residuals}, we are also able to bound $e_4$ as follow
\bequ\label{eq:CH3:e_4}
\Prob^{ \otimes n} \left\{ {e_4  > 0} \right\} \leqslant
    \sum\limits_{\textsf{I} \in \Lambda _{\MyJ^{\star} } } {\Prob^{ \otimes n}
    \left( {\{ \nu _\textsf{I} (\bz) \geqslant \lambda _n \}  \cap
    \{ \nu _\textsf{I}  \leqslant 2^{ - 1} \lambda _n \} } \right)}  \leqslant
    C'\,n^{ - \left[g(\kappa ) - \tfrac{1}{\gamma }\right]}.
\eequ
Combining Equations (\ref{eq:CH3:e_2_A}), (\ref{eq:CH3:e_2_B}), (\ref{eq:CH3:e_2_C}), and (\ref{eq:CH3:e_4}), we complete the proof of Equation (\ref{eq:CH3:e_2&e_4_goal}). In fact, given $\beta$ and $\gamma \geqslant \tfrac{1}{2}$, we can find $\kappa$ such that the theorem holds.

\subsubsection{Proof of Lemma \ref{lem:CH3:residuals}}

Lets starting noticing that, for each $\eta > 0$,
\[
\{ \nu _\textsf{I} (\bz) \leqslant \eta \}  \cap
    \{ \nu _\textsf{I}  \geqslant 2\eta \}  \subseteq
    \left\{ {|\nu_\textsf{I} (\bz) - \nu _\textsf{I} | \geqslant \eta } \right\},
\]
hence
\[
\Prob^{ \otimes n} \left( {\{ \nu _\textsf{I} (\bz) \leqslant \eta \}  \cap
    \{ \nu _\textsf{I} \geqslant 2 \eta \} } \right) \leqslant
    \Prob^{ \otimes n} \left\{ {|\nu _\textsf{I} (\bz) - \nu _\textsf{I} |
    \geqslant \eta } \right\}.
\]
In addition
\be
|\nu_\textsf{I}(\bz) - \nu_\textsf{I} | &=&
    \left| {\sqrt {\sum\nolimits_{\textsf{J} \in
    \mathcal{C}(\textsf{I})} {d_\textsf{I}^2(\bz)} }  -
    \sqrt {\sum\nolimits_{\textsf{J} \in \mathcal{C}(\textsf{I})} {d_\textsf{I}^2}}}\;\right| =
    \left| {\left\| {\mathbf{d}(\bz)} \right\|_2 - \left\| \mathbf{d} \right\|_2 } \right|
    \leqslant \hfill \\
&\leqslant&
    \left\| {\mathbf{d}(\bz) - \mathbf{d}} \right\|_2 \mathop  = \limits^{\text{dyadic}}
    \sqrt{
        \big[d_{\textsf{I}^+}(\bz) - d_{\textsf{I}^+}\big]^2 +
        \big[d_{\textsf{I}^-} (\bz) - d_{\textsf{I}^-}\big]^2
        },
\ee
where $\textsf{I}^+$ and $\textsf{I}^-$ denote respectively the left and right child of $\textsf{I}$. So
\be
\left\{ {|\nu _\textsf{I} (\bz) - \nu _\textsf{I} | \geqslant \eta } \right\}
    \Leftrightarrow \left\{ {|\nu_\textsf{I}(\bz) - \nu_\textsf{I} |^2 \geqslant\eta^2}\right\}
&\Leftarrow&
\bfun
   {[d_{I^ +  } (z) - d_{I^ +  } ]^2  \geqslant \tfrac{{\eta ^2 }}{2}}  \\
   {[d_{I^ -  } (z) - d_{I^ -  } ]^2  \geqslant \tfrac{{\eta ^2 }}{2}}  \\
\efun
\hfill \\
&\Leftrightarrow&
\bfun
   {|d_{I^ +  } (z) - d_{I^ +  } | \geqslant \tfrac{\eta }{{\sqrt 2 }}}  \\
   {|d_{I^ -  } (z) - d_{I^ -  } | \geqslant \tfrac{\eta }{{\sqrt 2 }}}  \\
\efun
.
\ee
Therefore
\[
\Prob^{ \otimes n} \left( {\{ \nu _\textsf{I}(\bz) \leqslant \lambda _n \}
    \cap \{ \nu _\textsf{I} \geqslant 2\lambda _n \} } \right) \leqslant
\Prob^{ \otimes n} \!\! \left( {| \Delta_{\textsf{I}^+}(\bz) | \geqslant
    \!\!\tfrac{\eta }{{\sqrt 2 }}} \right) +
\Prob^{ \otimes n} \!\! \left( {| \Delta_{\textsf{I}^-}(\bz) | \geqslant
    \!\!\tfrac{\eta }{{\sqrt 2 }}} \right).
\]
If we now take $\eta = \kappa \sqrt {\tfrac{{\log (n)}}{n}}$, by applying the Bernstein's inequality as in Section \ref{subsec:CH3:ProofUniform}, for $\textsf{J} \in \{\textsf{I}^+,\textsf{I}^-\}$ we obtain\footnote{Compare with the proof of Proposition 3 in \cite{KP04}.}
\be
\Prob^{ \otimes n} \left( {|\Delta _{\textsf{J}}(\bz)| \geqslant
    \tfrac{\kappa }{{\sqrt 2 }}\sqrt {\tfrac{{\log (n)}}{n}} } \right)
\!\!\!\!\! &\leqslant& \!\!\!\!\!
    2 \exp \left\{- {\frac{{3\kappa ^2 \log (n)}}{{2C'\big[3 + 2^{(j + 1)/2} 2^{ - 1/2} \kappa
        \sqrt {\tfrac{{\log (n)}}{n}} \big]}}} \right\} \leqslant \hfill \\
\!\!\!\!\! &\leqslant& \!\!\!\!\!
    2\exp \left\{- {\frac{{3\kappa ^2 \log (n)}}{{8C'\big[3 + 2^{j/2} \kappa
        \sqrt {\tfrac{{\log (n)}}{n}} \big]}}} \right\}.
\ee
Now, by hypothesis, we know that
\[
2^j  \leqslant \lambda _n^{ - 1/\gamma }  = \left[ {\tfrac{1}{\kappa }
    \sqrt {\tfrac{n}{{\log (n)}}} } \right]^{\tfrac{1}{\gamma }}
\Rightarrow
2^{j/2}  \leqslant \left[ {\tfrac{1}{\kappa }
    \sqrt {\tfrac{n}{{\log (n)}}} } \right]^{\tfrac{1}{{2\gamma }}},
\]
therefore
\[
\frac{{2^{j/2} }}{{\tfrac{1}{\kappa }\sqrt {\tfrac{n}{{\log (n)}}} }} \leqslant
    \left[ {\tfrac{1}{\kappa }\sqrt {\tfrac{n}{{\log (n)}}} } \right]^{
        - \left( {1 - \tfrac{1}{{2\gamma }}} \right)}  =
    \left[ {\kappa \sqrt {\tfrac{{\log (n)}}{n}} } \right]^{1 - \tfrac{1}{{2\gamma }}}
    \mathop  \leqslant \limits^{\gamma  \geqslant \tfrac{1}{2}} \kappa ^{1 - \tfrac{1}{{2\gamma}}},
\]
hence
\be
\Prob^{ \otimes n} \left( {|\Delta _\textsf{J}(\bz)| \geqslant
    \tfrac{\kappa }{{\sqrt 2 }}\sqrt {\tfrac{{\log (n)}}{n}} } \right) &\leqslant&
    2\exp \left\{ { - \frac{{3\kappa ^2 \log (n)}}{{8C'\big( {3 + \kappa ^{1
    - \tfrac{1}{{2\gamma }}} } \big)}}} \right\} = \hfill \\
    &=& 2\exp \left\{ { - g(\bz)\,\log (n)} \right\} = 2\,n^{ - g(\bz)},
\ee
with
\[
g(\kappa ) = \frac{{3\kappa ^2 }}{{8C'\biggl( {3 + \kappa ^{1 - \tfrac{1}{{2\gamma }}} }\biggr)}}.
\]
So finally
\[
\Prob^{ \otimes n} \left( {\{ \nu _\textsf{I} (\bz) \leqslant \lambda _n \}
    \cap \{ \nu _\textsf{I} \geqslant 2\lambda _n \} } \right) \leqslant 4\,n^{ - g(\bz)} .
\]
Now, lets evaluate the other term in a similar manner, starting from
\be
\Prob^{ \otimes n} \left( {\{ \nu _\textsf{I} (\bz) \geqslant \lambda _n \}  \cap
    \{ \nu_\textsf{I} \leqslant 2^{ - 1} \lambda _n \} } \right) &\leqslant&
    \Prob^{ \otimes n} \left\{ {|\nu_\textsf{I}(\bz) - \nu _\textsf{I} | \geqslant
    2^{ - 1} \lambda _n } \right\} \leqslant \hfill \\
&\leqslant&
    \sum\limits_{\textsf{J} \in \{ \textsf{I}^ + ,\textsf{I}^ -  \} }
    {\Prob^{ \otimes n} \left( {|\Delta _\textsf{J} (\bz)| \geqslant
    \tfrac{{\lambda _n }}{{2\sqrt 2 }}} \right)}.
\ee
By the same arguments adopted before, we see that, for each $\textsf{J} \in \{\textsf{I}^+,\textsf{I}^-\}$,
\[
\Prob^{ \otimes n} \left( {|\Delta _\textsf{J} (\bz)| \geqslant \tfrac{{\lambda _n }}
{{2\sqrt 2 }}} \right) \leqslant 2\,n^{ - g(\bz)},
\]
and consequently
\[
\Prob^{ \otimes n} \left( {\{ \nu _\textsf{I} (\bz) \geqslant \lambda _n \}  \cap \{ \nu _\textsf{I} \leqslant 2^{ - 1} \lambda _n \} } \right) \leqslant 4\,n^{ - g(\bz)} .
\]

\bibliographystyle{plain}
\def\cprime{$'$}

\end{document}